\documentclass{elsart}
\usepackage{amssymb,amsfonts}
\newcommand{\trace}{\mathop{\rm Tr}\nolimits}

\newcommand{\vect}[1]{\mathbf{#1}}

\newcommand{\re}{\mathop{\rm Re}\nolimits}

\newcommand{\twomat}[4]{\left(\begin{array}{cc}#1&#2\\#3&#4\end{array}\right)}

\newcommand{\cP}{{\mathcal P}} 

\newcommand{\C}{{\mathbb{C}}}
\newcommand{\R}{{\mathbb{R}}}

\newcommand{\be}{\begin{equation}}
\newcommand{\ee}{\end{equation}}
\newcommand{\bea}{\begin{eqnarray}}
\newcommand{\eea}{\end{eqnarray}}
\newcommand{\beas}{\begin{eqnarray*}}
\newcommand{\eeas}{\end{eqnarray*}}
\newtheorem{definition}{Definition}
\newtheorem{theorem}{Theorem}
\newtheorem{lemma}{Lemma}

\newcount\minute
\newcount\hour
\def\currenttime{%
    \minute\time
    \hour\minute
    \divide\hour60
    \the\hour:\multiply\hour60\advance\minute-\hour\the\minute}
\begin{document}
\begin{frontmatter}
\title{A Determinantal Inequality for the Geometric Mean with an Application in Diffusion Tensor Imaging}
\author{Koenraad M.R.~Audenaert}
\address{
Department of Mathematics,
Royal Holloway University of London, \\
Egham TW20 0EX, United Kingdom \\[1mm]
Department of Physics and Astronomy, Ghent University, \\
S9, Krijgslaan 281, B-9000 Ghent, Belgium}
\ead{koenraad.audenaert@rhul.ac.uk}
\date{\today, \currenttime}
\begin{abstract}
We prove that for positive semidefinite matrices $A$ and $B$
the following determinantal inequality holds:
\[
\det(I+A\#B)\le \det(I+A^{1/2}B^{1/2}),
\]
where $A\#B$ is the geometric mean of $A$ and $B$.
We apply this inequality to the study of interpolation methods in diffusion tensor imaging.
\end{abstract}

\end{frontmatter}
\section{Introduction}
The topic of this paper has arisen in the study of interpolation methods
for image processing in diffusion tensor imaging (DTI). DTI is an imaging method used in medical magnetic resonance imaging (MRI)
whereby the data to be imaged usually consists of a field of $3\times 3$ statistical covariance matrices $D(\vect r)\in \cP$,
where the points $\vect r$ lie on a 3-dimensional grid.
In order to improve the visual quality of the image it is necessary to interpolate and/or extrapolate between neighbouring $D$-values.
As there is no unique way to do so, it is important to have a mathematical framework within which to describe the various methods
as well as their measures of goodness.
Such a framework has recently been discussed in \cite{dryden,alexey}, and the present work grew out of this.

The interpolation/extrapolation problem is easily formulated as follows. In this context, all covariance matrices are real, symmetric, positive
semidefinite $3\times 3$ matrices.
Let $D_1$ and $D_2$ be two covariance matrices. Construct a path $p\mapsto D(p)$, where $p\in[0,1]$ for interpolation or $p<0$ or $p>1$ for
extrapolation, such that $D(0)=D_1$, $D(1)=D_2$ and
$D(p)\ge0$, for all $p$ in the interval of interest, and whereby certain quality criteria have to be satisfied.
Without going in too much detail, the quality criterion used in \cite{alexey} is based on the cube root of the determinant of $D(p)$ as this is one
of the quantities that provides structural details of the sample being imaged.
It is argued by some that interpolated/extrapolated values of this determinant should not be ``too large'', as this
might lead to so-called ``swelling'' of certain features
in the reconstructed image.

One further requirement on the path $D(p)$ is that it be the shortest path between $D_1$ and $D_2$, in a sense to be specified.
This is easiest to satisfy by defining a suitable metric $d(\cdot,\cdot)$ on $\cP$ and let the path be a geodesic path.
The requirement that all $D(p)$ on the path be positive semidefinite can then be enforced by choosing a metric specific to the curved space $\cP$;
this excludes the Euclidean metric $d(A,B)=||A-B||_2$ for example as its geodesics do not stay within $\cP$.
Still, many choices remain. Probably the most studied $\cP$-metric in matrix analysis is the Riemannian metric
$d_R = ||\log(D_1^{-1/2}D_2D_1^{-1/2})||_2$ \cite[Chapter 6]{bhatiapdm}.
Accordingly, this metric has been given due consideration in DTI as well \cite{fletcher}.
However, one of its drawbacks is that it is inordinately sensitive to very small eigenvalues of
the covariance matrices. In particular, all rank-deficient matrices  are infinitely far apart in the Riemannian
sense, no matter how close they are in the Euclidean sense. For this and other reasons, other metrics besides the
Riemannian one are being studied for DTI.

The metrics studied in the present work are the so-called ``Euclidean root metric'', $d_H$, and the ``Procrustes size-and-shape metric'' $d_S$;
we will not explain this nomenclature here.
Both metrics are based on a reparameterisation of the covariance matrices to enforce positive semidefiniteness, in that they
give the Euclidean distance between certain square roots of the covariance matrices, $||S_1-S_2||_2$.
Here, a square root of a positive semidefinite matrix $D$ is any matrix $S$ for which $D=S^*S$.
Taking for example the Cholesky decomposition for $S$ ($S$ being upper triangular with positive diagonal entries) yields the Cholesky metric $d_C$.
One gets the Euclidean root metric $d_H$ by taking positive square roots; i.e.\ $S_i=Q_i:=D_i^{1/2}$.

The Procrustes metric is the minimal one in the sense that it minimises $||S_1-S_2||_2$ over all possible choices of square roots.
That is,
\[
d_S(D_1,D_2) = \min_R\{||Q_1-RQ_2||_2: R \mbox{ unitary}\}.
\]
In other words, we look for that square root of $D_2$ that is closest to the positive square root of $D_1$.
This minimisation can easily be done analytically, since $||Q_1-RQ_2||_2^2 = \trace (Q_1^2+Q_2^2)-2\re\trace RQ_2Q_1$.
Hence, the optimal $R$ in the above minimisation is the unitary matrix for which $RQ_2Q_1 = |Q_2Q_1|$, where $|X|=\sqrt{X^*X}$ denotes
the modulus of the matrix $X$.
Thus the optimal $R$ is given by $R=U^*$
where $U$ is the unitary factor in the polar decomposition of $Q_2Q_1=U|Q_2Q_1|$.

The geodesics induced by these metrics are obtained by considering linear paths in the square root space, and taking the square to map them back to
$\cP$.
That is, for the Euclidean root metric,
\[
D_H(p) = |pQ_1+(1-p)Q_2|^2
\]
and for the Procrustes metric
\[
D_S(p) = |pQ_1+(1-p)U^*Q_2|^2.
\]

The question that we wanted to answer about these metrics is how the determinants of the interpolated values of the various $D(p)$ behave.
Our main result shows that for $0\le p\le 1$ the Procrustes path always produces determinants that are smaller or equal to those of the Euclidean root path,
i.e.\ the Procrustes interpolation is less prone to swelling. For the extrapolation ($p<0$ or $p>1$) numerical calculations confirmed
our intuition that there can be no guaranteed ordering between the determinants of the two paths.

\section{Main Result}

The following theorem answers the question posed in the introduction; in fact it is more general than required as it holds for all
complex positive semidefinite matrices, of any dimension.
\begin{theorem}\label{th:main}
Let $Q_1$ and $Q_2$ be $n\times n$ positive semidefinite matrices,
Let $Q_2Q_1$ have polar decomposition $Q_2Q_1=U|Q_2Q_1|$.
Then $\det(Q_1+U^*Q_2)$ is real and non-negative and
\[
\det(Q_1+U^*Q_2) \le \det(Q_1+Q_2).
\]
\end{theorem}
Here we have absorbed the interpolation parameter $p$ in $Q_1$ and $1-p$ in $Q_2$. For $0\le p\le 1$ this does not change the signature of
$Q_1$ or $Q_2$, nor does it change the unitary factor $U$.

That $\det(Q_1+U^*Q_2)$ is real and non-negative is easy to show.
Indeed, we have $U^*Q_2Q_1 = |Q_2Q_1|$, so that $\det(Q_1+U^*Q_2)\det(Q_1)=\det(Q_1^2+|Q_2Q_1|)>=0$. Dividing by the positive number $\det(Q_1)$ shows
that $\det(Q_1+U^*Q_2)\ge0$.

The proof of the inequality of this theorem relies on a number of concepts from matrix analysis, which we introduce first in some detail
(for the benefit of those readers working in diffusion tensor imaging).
In Section \ref{sec:majo} we consider the related concepts of majorisation and $\log$-majorisation, presenting their definitions
and their most important properties.
In Section \ref{sec:geomean} we consider the matrix generalisation of the geometric mean.
Finally, in Section \ref{sec:proof} we present the proof of the theorem.

\section{Majorisation\label{sec:majo}}
In this section we consider vectors $\vect x=(x_1,\ldots,x_n)$ and $\vect y=(y_1,\ldots,y_n)$ in $\R^n$.
We denote by $\vect x^\downarrow$ the vector consisting of the elements of $\vect x$, sorted in non-ascending order.
Thus, $x^\downarrow_k$ is the $k$-th largest element of $\vect x$.

We will now introduce several relations between $\vect x$ and $\vect y$ that come under the heading of \emph{majorisation}.
The standard work about the theory and applications of majorisation is undoubtedly \cite{MO}, to which we refer for more details.
A more concise treatment can be found in \cite[Chapter II]{bhatia}.

We say that $\vect x$ \emph{weakly majorises} $\vect y$, denoted $\vect y\prec_w \vect x$,
if and only if, for all $k$, the sum of the $k$ largest
elements of $\vect x$ dominates the sum of the $k$ largest elements of $y$:
\[
\vect y\prec_w \vect x \Longleftrightarrow \sum_{i=1}^k y^\downarrow_i \le \sum_{i=1}^k x^\downarrow_i, \quad k=1,2,\ldots,n.
\]
If in addition the sum of all elements of $\vect x$ equals the sum of all elements of $\vect y$, then we say
that $\vect x$ \emph{majorises} $\vect y$:
\[
\vect y\prec \vect x \Longleftrightarrow
\left\{
\begin{array}{l}
\sum_{i=1}^k y^\downarrow_i \le \sum_{i=1}^k x^\downarrow_i, \quad k=1,2,\ldots,n-1;\\
\sum_{i=1}^n y^\downarrow_i = \sum_{i=1}^n x^\downarrow_i.
\end{array}
\right.
\]

A closely related concept is log-majorisation, which concerns vectors in $\R_+^n$.
We say that $\vect x$ \emph{weakly log-majorises} $\vect y$, denoted $\vect y\prec_{w,\log}\vect x$
if and only if $\log\vect y\prec_w\log\vect x$.
Thus,
\[
\vect y\prec_{w,\log}\vect x \Longleftrightarrow \prod_{i=1}^k y^\downarrow_i \le \prod_{i=1}^k x^\downarrow_i, \quad k=1,2,\ldots,n.
\]
Likewise, $\vect x$ \emph{log-majorises} $\vect y$, denoted $\vect y\prec_{\log}\vect x$
if and only if $\log\vect y\prec\log\vect x$.

Next we discuss which functions preserve the (weak) majorisation ordering.
A function $\Phi:\R^n\to\R^m$ is called \emph{strongly isotone} if and only if it preserves weak majorisation:
\[
\vect y \prec_w \vect x \Rightarrow \Phi(\vect y) \prec_w \Phi(\vect x).
\]
A function is called \emph{isotone} if and only if
\[
\vect y \prec \vect x \Rightarrow \Phi(\vect y) \prec_w \Phi(\vect x).
\]
We will need the following characterisation of isotony in the case that $m=1$ \cite[Theorem II.3.14]{bhatia}.
\begin{lemma}[Schur]
A differentiable function $\Phi:\R^n\to\R$ is isotone if and only if it satisfies
\begin{enumerate}
\item $\Phi$ is permutation invariant,
\item for all $\vect x\in \R^n$ and for all $i,j$:
\[
(x_i-x_j)\left(\frac{\partial\Phi}{\partial x_i}(\vect x) - \frac{\partial\Phi}{\partial x_j}(\vect x)\right)\ge0.
\]
\end{enumerate}
\end{lemma}
\section{Geometric Mean\label{sec:geomean}}

Recall that the geometric mean of two positive real numbers $x$ and $y$ is given by $\sqrt{xy}$.
As is the case with all means, when one wishes to generalise the geometric mean to two positive semidefinite matrices $A$ and $B$,
one is faced with the usual problem of non-commutativity: there exists an infinite number of expressions involving $A$ and $B$
that reduce to $\sqrt{AB}$ when $A$ and $B$ commute. For example, $\sqrt{AB}$ and $\sqrt{A}\sqrt{B}$ are in general
different matrices, but both are equal in the commutative case.

To resolve this problem, one has to impose a number of conditions in order to obtain
a unique generalisation.
Kubo and Ando \cite{KA} have developed a very nice theory of matrix means that does exactly that.
It is now standard to define the matrix geometric mean as follows:
\begin{definition}
Let $A$ and $B$ be $n\times n$ positive definite matrices. Then their geometric mean, denoted $A\#B$, is defined as
\[
A\#B := A^{1/2} (A^{-1/2} B A^{-1/2})^{1/2} A^{1/2}.
\]
When $A$ and/or $B$ are rank-deficient, the geometric mean is defined via a limiting procedure.
\end{definition}
While this is not obvious from the definition, the geometric mean is symmetric in its arguments; that is, $A\#B=B\#A$.

Clearly, for $a,b>0$ we have $(aA)\#(bB)=\sqrt{ab}(A\#B)$.

One can show that $A\#B$ emerges as the solution of the following optimisation problem:
the set of positive semidefinite matrices $X$ for which the block matrix
\[
\twomat{A}{X}{X}{B}
\]
is itself positive semidefinite, has a maximum in the positive semidefinite ordering, and this maximum is exactly $A\#B$.

In the present work, we will need two more properties of the geometric mean:
\begin{lemma}[Monotonicity]
Let $0\le A$ and $0\le B_1\le B_2$. Then
\[
A\#B_1 \le A\#B_2.
\]
\end{lemma}

\begin{lemma}\label{lem:hiai}
Let $A,B\ge 0$ and $r\ge 1$.
If $A\#B \le 1$, then $A^r\#B^r\le 1$.
\end{lemma}
The proof of this last statement can be found in \cite{hiai}, as part of the proof of its Theorem 4.6.9.

\section{Proof of the theorem.\label{sec:proof}}

We start with a lemma relating the eigenvalues of the matrix geometric mean to the eigenvalues of
what could be designated as a ``na{\"\i}ve'' matrix geometric mean.
This is actually a special case of a more general inequality due to Matharu and Aujla \cite[Theorem 2.10]{matharu}, but we provide
a stand-alone proof for the benefit of the reader.
Henceforth we will denote the vector of eigenvalues of a matrix $X$ by $\lambda(X)$.
\begin{lemma}\label{lem:mine}
Let $A,B\ge0$.
Then
\[
\lambda(A\#B)\prec_{\log} \lambda(\sqrt{A}\;\sqrt{B}).
\]
\end{lemma}
\textit{Proof.}
Throughout the proof we assume that $A$ is invertible.
The general case follows from continuity considerations.

We first show that the inequality holds for the largest eigenvalue $\lambda_1$; that is:
\be
\lambda_1(A\#B)\le \lambda_1(\sqrt{A}\;\sqrt{B}).\label{eq:l1}
\ee
Let $a=\lambda_1(\sqrt{A}\;\sqrt{B})$. Thus $B^{1/2}\le a A^{-1/2}$.
By monotonicity of the geometric mean, we then have
$A^{1/2} \# B^{1/2} \le \sqrt{a}\;(A^{1/2} \# A^{-1/2}) = \sqrt{a}$.
Using Lemma \ref{lem:hiai} with $r=2$, we obtain $A\#B \le a$.
This says that $\lambda_1(A\#B)\le a$, as required.

To prove the statement of the lemma, we use the so-called ``Weyl-trick''.
Let $A^{\wedge k}$ be the $k$-th antisymmetric tensor power of $A$; this is the restriction
of the $k$-th tensor power of $A$, $A^{\otimes k}$, to the totally antisymmetric subspace of $(\C^n)^{\otimes k}$.
The Weyl-trick exploits two facts about these powers. Firstly, for $A\ge0$ the largest eigenvalue of $A^{\wedge k}$ is given by
\[
\lambda_1(A^{\wedge k}) = \lambda_1(A)\lambda_2(A)\cdots\lambda_k(A).
\]
Secondly, any expression involving products and/or fractional matrix powers ``commutes'' with taking the $k$-th
antisymmetric tensor power.
In particular, $(A\#B)^{\wedge k} = A^{\wedge k}\#B^{\wedge k}$.

Thus, exploiting (\ref{eq:l1}), we get
\beas
\prod_{i=1}^k \lambda_i(A\#B) &=& \lambda_1((A\#B)^{\wedge k}) = \lambda_1(A^{\wedge k} \# B^{\wedge k}) \\
&\le& \lambda_1(\sqrt{A^{\wedge k}}\sqrt{B^{\wedge k}}) = \lambda_1((\sqrt{A}\;\sqrt{B})^{\wedge k})
= \prod_{i=1}^k \lambda_i(\sqrt{A}\;\sqrt{B}).
\eeas
This already shows that we have a weak log-majorisation: $\lambda(A\#B)\prec_{w,\log} \lambda(\sqrt{A}\;\sqrt{B})$.

Because of the Cauchy-Binet theorem, we also have
\[
\det(A\#B) = \sqrt{\det(A)\det(B)} = \det(\sqrt{A}\sqrt{B}),
\]
which says that
$\prod_{i=1}^n \lambda_i(A\#B) = \prod_{i=1}^n \lambda_i(\sqrt{A}\;\sqrt{B})$.
Thus, the weak log-majorisation relation can be strengthened to a log-majorisation.
\qed

\begin{theorem}\label{th:detgm}
Let $A,B\ge0$. Then
\be
\det(I+A\#B) \le \det(I + A^{1/2}B^{1/2}).
\ee
\end{theorem}
\textit{Proof.}
Note that, for any $X\ge 0$, we have $\log\det X = \trace\log X=\sum_{i=1}^n \log\lambda_i(X)$.
So we equivalently need to show
\be
\trace\log(I+A\#B) \le \trace\log(I + A^{1/2}B^{1/2}).\label{eq:rephrase}
\ee

By Lemma \ref{lem:mine},
\[
\log\lambda(A\#B)\prec \log\lambda(A^{1/2}B^{1/2}).
\]
This implies (\ref{eq:rephrase}) if we can show that the function
$\Phi(\vect x):=\sum_{i=1}^n \log(1+\exp(x_i))$ is isotone.

Clearly, this function is permutation symmetric.
Furthermore,
\[
\partial\Phi/\partial x_i = \exp(x_i)/(1+\exp(x_i)) = 1/(1+\exp(-x_i)),
\]
which is monotonously increasing in $x_i$.
Hence, the second condition for isotony is also satisfied.
\qed

\textit{Proof of Theorem \ref{th:main}.}\\
Let us consider the matrices $Q_1$, $Q_2$ and $U$ given in the statement of the theorem.
Thus $U^*Q_2Q_1 = |Q_2Q_1|$.
We will assume that $Q_1$ is invertible, so $\det(Q_1)\neq 0$.
The general case follows from continuity of the determinant.

Because of the Cauchy-Binet theorem, the statement of the theorem is equivalent with
\[
\det((Q_1+U^*Q_2)Q_1) \le \det(Q_1(Q_1+Q_2)),
\]
which becomes
\[
\det(Q_1^2+|Q_2Q_1|) \le \det(Q_1^2+Q_1Q_2).
\]
Applying the Cauchy-Binet theorem a second time we can divide out $Q_1^2$ from each side in a well-chosen way, and get
the equivalent statement
\[
\det(I+Q_1^{-1}(Q_1Q_2^2Q_1)^{1/2}Q_1^{-1}) \le \det(I+Q_1^{-1/2}Q_2Q_1^{-1/2}).
\]
Using the substitutions $A=Q_1^{-2}$ and $B=Q_2^2$, this can be rewritten as
\[
\det(I+A^{1/2}(A^{-1/2}BA^{-1/2})^{1/2}A^{1/2}) \le \det(I+A^{1/4}B^{1/2}A^{1/4}),
\]
or, in terms of the geometric mean,
\[
\det(I+A\#B) \le \det(I + A^{1/2}B^{1/2}).
\]
The validity of this inequality is just Theorem \ref{th:detgm}.
\qed

\section*{Acknowledgments}
We acknowledge support by an Odysseus grant from the Flemish FWO.
Many thanks to A.~Koloydenko for suggesting this problem.


\end{document}